\input{style/arxiv-general.cfg}
\documentclass[MSNbibl,number,citesort,seceqn,dvips]{arxbj}
\makeatletter
   \@ifpackageloaded{graphicx}{}{\usepackage{graphicx}}
\makeatother
\usepackage{dcolumn}
\usepackage{upgreek}
\usepackage{graphicx}

\aid{0}
\volume{21}
\issue{2}
\pubyear{2015}
\firstpage{909}
\lastpage{929}
\doi{10.3150/13-BEJ592} 

\makeatletter
\newcolumntype{d}[1]{D{.}{.}{#1}}
\newremark{remark}{Remark}[section]
\newproclaim{definition}{Definition}[section]
\newproclaim{assumption}{Assumption}[section]
\newtheorem{theorem}{Theorem}[section]
\newtheorem{proposition}{Proposition}[section]
\newcommand{\field}[1]{\mathbb{#1}}
\newcommand{\R}{\field{R}}
\newcommand{\toprule}{\hline}
\newcommand{\midrule}{\hline}
\makeatother

\begin{document}
\begin{frontmatter}

\title{Two sample inference for the second-order property of
temporally dependent functional data}
\runtitle{Two sample inference}

\begin{aug}
\author[A]{\inits{X.}\fnms{Xianyang} \snm{Zhang}\corref{}\thanksref{A}\ead[label=e1]{zhangxiany@missouri.edu}} 
\and
\author[B]{\inits{X.}\fnms{Xiaofeng} \snm{Shao}\thanksref{B}\ead[label=e2]{xshao@illinois.edu}}
\address[A]{Department of Statistics, University of Missouri-Columbia,
Columbia, MO, 65211, USA.\\ \printead{e1}}
\address[B]{Department of Statistics, University of Illinois at
Urbana-Champaign, Champaign, IL, 61820, USA. \printead{e2}}
\end{aug}

\received{\smonth{12} \syear{2012}}
\revised{\smonth{12} \syear{2013}}

%
\begin{abstract}
Motivated by the need to statistically quantify the difference
between two spatio-temporal datasets that arise in climate
downscaling studies, we propose new tests to detect the differences
of the covariance operators and their associated characteristics of
two functional time series. Our two sample tests are constructed on
the basis of functional principal component analysis and
self-normalization, the latter of which is a new studentization
technique recently developed for the inference of a univariate time
series. Compared to the existing tests, our SN-based tests allow for
weak dependence within each sample and it is robust to the
dependence between the two samples in the case of equal sample
sizes. Asymptotic properties of the SN-based test statistics are
derived under both the null and local alternatives. Through
extensive simulations, our SN-based tests are shown to outperform
existing alternatives in size and their powers are found to be
respectable. The tests are then applied to the gridded climate model
outputs and interpolated observations to detect the difference in
their spatial dynamics.
\end{abstract}

%

\end{frontmatter}

\section{Introduction}

Functional data analysis (FDA) which deals with the analysis of
curves and surfaces has received considerable attention in the
statistical literature during the last decade (Ramsay and Silverman
\cite{r22,r23} and Ferraty and Vieu \cite{r8}).
This paper falls into a
sub-field of functional data analysis: inference for temporally
dependent functional data. Specifically, we focus on testing the
equality of the second-order structures (e.g., the covariance
operators and their associated eigenvalues and eigenfunctions) of
two temporally dependent functional sequences. Our work is partially
motivated by our ongoing collaboration with atmospheric scientists
on the development and assessment of high-resolution climate
projections through statistical downscaling. Climate change is one
of the most urgent problems facing the world this century. To study
climate change, scientists have relied primarily on climate
projections from global/regional climate models, which are numerical
models that involve systems of differential equations and produce
outputs at a prespecified grid. As numerical model outputs are
widely used in situations where real observations are not available,
it is an important but still open question whether the numerical
model outputs are able to mimic/capture the spatial and temporal
dynamics of the real observations. To partly answer this question,
we view the spatio-temporal model outputs and real observations as
realizations from two temporally dependent functional time series
defined on the two-dimensional space and test the equality of their
second-order structures which reflects their spatial
dynamics/dependence.

Two sample inference for functional data has been investigated by a
few researchers. Fan and Lin \cite{r7}, Cuevas \textit{et al.} \cite{r6} and
Horv\'{a}th \textit{et al.} \cite{r12} developed the tests for the equality of
mean functions. Benko \textit{et al.} \cite{r2}, Panaretos \textit{et al.}
\cite{r19}, Fremdt \textit{et al.} \cite{r9}, and Kraus and Panaretos
\cite{r15} proposed tests for the
equality of the second-order structures. All the above-mentioned
works assumed the independence between the two samples and/or
independence within each sample. However, the assumption of
independence within the sample is often too strong to be realistic
in many applications, especially if data are collected sequentially
over time. For example, the independence assumption is questionable
for the climate projection data considered in this paper, as the
model outputs and real station observations are simulated or
collected over time and temporal dependence is expected. Furthermore
the dependence between numerical model outputs and station
observations is likely because the numerical models are designed to
mimic the dynamics of real observations. See Section~\ref{data} for
empirical evidence of their dependence. In this paper, we develop
new tests that are able to accommodate weak dependence between and
within two samples. Our tests are constructed on the basis of
functional principal component analysis (FPCA) and the recently
developed self-normalization (SN) method (Shao \cite{r24}), the latter
of which is a new studentization technique for the inference of a
univariate time series.

FPCA attempts to find the dominant modes of variation around an
overall trend function and has been proved a key technique in the
context of FDA. The use of FPCA in the inference of temporally
dependent functional data can be found in Gabrys and Kokoszka \cite
{r10}, H\"{o}rmann and Kokoszka \cite{r11}, Horv\'{a}th \textit{et al.}
\cite{r12}
among others. To account for the dependence, the existing inference
procedure requires a consistent estimator of the long run variance
(LRV) matrix of the principal component scores or consistent
estimator of the LRV operator. However, there is a bandwidth
parameter involved in the LRV estimation and its selection has not
been addressed in the functional setting. The same issue appears
when one considers the block bootstrapping and subsampling schemes
(Lahiri \cite{r16} and Politis \textit{et al.} \cite{r21}),
since these techniques also
require the selection of a smoothing parameter, such as the block
length in the moving block bootstrap, and the window width in the
subsampling method (see, e.g., Politis and Romano \cite{r20} and
McMurry and Politis \cite{r18}). Since the finite sample performance
can be
sensitive to the choice of these tuning parameters and the bandwidth
choice can involve certain degree of arbitrariness, it is desirable
to use inference methods that are free of bandwidth parameters. To
this end, we build on the bandwidth-free SN method (Shao
\cite{r24})
recently developed in the univariate time series setup, and propose
SN-based tests in the functional setting by using recursive
estimates obtained from functional data samples.


In time series analysis, the inference of a parameter using normal
approximation typically requires consistent estimation of its
asymptotic variance. The main difficulty with this approach (and
other block-based resampling methods) is the sensitivity of the
finite sample performance with respect to the bandwidth parameter,
which is often difficult to choose in practice without any
parametric assumptions. As a useful alternative, the
self-normalized approach uses an inconsistent bandwidth-free
estimate of asymptotic variance, which is proportional to asymptotic
variance, so the studentized quantity (statistic) is asymptotically
pivotal. Extending the early idea of Lobato \cite{r17},
Shao \cite{r24}
proposed a very general kind of self-normalizers that are functions
of recursive estimates and showed the theoretical validity for a
wide class of parameters of interest. The settings in the latter two
papers are however limited to univariate time series. The
generalization of the SN method from univariate to functional time
series was first done in Zhang \textit{et al.} \cite{r26} where the
focus was on
testing the structure stability of temporally dependent functional
data. Here we extend the SN method to test the equality of the
second-order properties of two functional time series, which is
rather different and new techniques and results are needed. To study
the asymptotic properties of the proposed test statistics, we
establish functional central limit theorems for the recursive
estimates of quantities associated with the second-order properties
of the functional time series which seems unexplored in the
literature and are thus of independent interest. Based on the
functional central limit theorem, we show that the SN-based test
statistics have pivotal limiting distributions under the null and
are consistent under the local alternatives. From a methodological
viewpoint, this seems to be the first time that the SN method is
extended to the two sample problem. Compared to most of the existing
methods which assumed the independence between the two samples
and/or independence within each sample, the SN method not only
allows for unknown dependence within each sample but also allows for
unknown dependence between the two samples when the sample sizes of
the two sequences are equal.

\section{Methodology}

We shall consider temporally dependent functional processes
$\{(X_i(t),Y_i(t)),t\in\mathcal{I}\}_{i=1}^{+\infty}$ defined on
some compact set $\mathcal{I}$ of the Euclidian space, where
$\mathcal{I}$ can be one-dimensional (e.g., a curve) or
multidimensional (e.g., a surface or manifold). For simplicity, we
consider the Hilbert space $\mathbb{H}$ of square integrable
functions with $\mathcal{I}=[0,1]$ (or $\mathcal{I}=[0,1]^2$). For
any functions $f,g\in\mathbb{H}$, the inner product between $f$ and
$g$ is defined as $\int_{\mathcal{I}}f(t)g(t)\,\mathrm{d}t$ and
$\Vert\cdot\Vert$
denotes the inner product induced norm. Assume the random elements
all come from the same probability space $(\Omega, \mathcal{A},
\mathcal{P})$. Let $L^{p}$ be the space of real valued random
variables with finite $L^{p}$ norm, that is, $(E|X|^p)^{1/p}<\infty$.
Further, we denote $L^p_{\mathbb{H}}$ the space of $\mathbb{H}$
valued random variables $X$ such that $(E\Vert X\Vert^p)^{1/p}<\infty$.

Given two sequences of temporally dependent functional observations,
$\{X_i(t)\}_{i=1}^{N_1}$ and $\{Y_i(t)\}_{i=1}^{N_2}$ defined on a
common region $\mathcal{I}$, we are interested in comparing their
second-order properties. Suppose that the functional time series are
second-order stationary. We assume that $E[X_i(t)]=E[Y_i(t)]=0$. The
result can be easily extended to the situation with nonzero mean
functions. Define $C_X=E[\langle X_i,\cdot\rangle X_i]$ and
$C_Y=E[\langle Y_i,\cdot\rangle Y_i]$ as the covariance operators of the two
sequences respectively. For the convenience of presentation, we
shall use the same notation for the covariance operator and the
associated covariance function. Denote by
$\{\phi_X^{j}\}_{j=1}^{\infty}$ and
$\{\lambda_{X}^{j}\}_{j=1}^{\infty}$ the eigenfunctions and
eigenvalues of $C_X$. Analogous quantities are
$\{\phi_Y^{j}\}_{j=1}^{\infty}$ and
$\{\lambda_{Y}^{j}\}_{j=1}^{\infty}$ for the second sample. Denote
by $|\mathbf{v}|$ the Euclidean norm of a vector
$\mathbf{v}\in\mathbb{R}^p$. Let $\operatorname{vech}(\cdot)$ be the
operator that stacks the columns below the diagonal of a symmetric
$m\times m$ matrix as a vector with $m(m+1)/2$ components. Let
$D[0,1]$ be the space of functions on $[0,1]$ which are
right-continuous and have left limits, endowed with the Skorokhod
topology (see Billingsley \cite{r3}). Weak convergence in
$D[0,1]$ or
more generally in the $\mathbb{R}^m$-valued function space $D^m[0,
1]$ is denoted by ``$\Rightarrow$'', where $m\in\mathbb{N}$ and
convergence in distribution is denoted by ``$\rightarrow^d$''. Define
$\lfloor a\rfloor$ the integer part of $a\in\mathbb{R}$, and
$\delta_{ij}=1$ if $i=j$ and $\delta_{ij}=0$ if $i\neq j$. In what
follows, we shall discuss the tests for comparing the three
quantities $C_X$, $\phi_X^{j}$ and $\lambda_X^{j}$ with $C_Y$,
$\phi_Y^{j}$ and $\lambda_Y^{j}$, respectively.

\subsection{Covariance operator}\label{test-cov-oper}

Consider the problem of testing the hypothesis $H_{1,0}\dvtx  C_X=C_Y$
versus the alternative $H_{1,a}\dvtx  C_X\neq C_Y$ (in the operator
norm sense) for two mean zero stationary functional time series
$\{X_i(t)\}_{i=1}^{N_1}$ and $\{Y_i(t)\}^{N_2}_{i=1}$. Let
$N=N_1+N_2$. Throughout the paper, we assume that
\[
N_1/N\rightarrow\gamma_1,\qquad N_2/N
\rightarrow\gamma_2, \qquad\mbox{as } \min(N_1,N_2)
\rightarrow+\infty,
\]
where $\gamma_1,\gamma_2\in(0,1)$ and $\gamma_1+\gamma_2=1$. Define
the one-dimensional operator
$\mathcal{X}_i=\langle X_i,\cdot\rangle X_i=X_i\otimes X_i$ and
$\mathcal{Y}_j=\langle Y_j,\cdot\rangle Y_j=Y_j\otimes Y_j$. Let $\hat{C}_{XY}$ be
the empirical covariance operator based on the pooled samples, that is,
%
%
\begin{equation}
\label{cov-poolsample} \hat{C}_{XY}=\frac{1}{N_1+N_2} \Biggl(\sum
^{N_1}_{i=1}\mathcal{X}_i +\sum
^{N_2}_{i=1}\mathcal{Y}_i \Biggr).
\end{equation}
Denote by $\{\hat{\lambda}_{XY}^{j}\}$ and $\{\hat{\phi}_{XY}^{j}\}$
the corresponding eigenvalues and eigenfunctions. The population
counterpart of $\hat{C}_{XY}$ is then given by
$\tilde{C}_{XY}=\gamma_1C_X+\gamma_2C_Y$ whose eigenvalues and
eigenfunctions are denoted by $\{\tilde{\lambda}^j\}$ and
$\{\tilde{\phi}^j\}$ respectively. Further let
$\hat{C}_{X,m}=\frac{1}{m}\sum^{m}_{i=1}\mathcal{X}_i$ be the sample
covariance operator based on the subsample $\{X_i(t)\}^{m}_{i=1}$
with $2 \leq m\leq N_1$. Define $\{\hat{\phi}_{X,m}^{j}\}^{m}_{j=1}$
and $\{\hat{\lambda}_{X,m}^{j}\}^{m}_{j=1}$ the eigenfunctions and
eigenvalues of $\hat{C}_{X,m}$ respectively, that is,
%
%
\begin{equation}
\label{eigenequation} \int_{\mathcal{I}}\hat{C}_{X,m}(t,s)\hat{
\phi}^{j}_{X,m}(s)\, \mathrm{d}s=\hat {\lambda}_{X,m}^j
\hat{\phi}_{X,m}^j(t),
\end{equation}
and
$\int_{\mathcal{I}}\hat{\phi}_{X,m}^i(t)\hat{\phi
}_{X,m}^j(t)\,\mathrm{d}t=\delta_{ij}$.
Similarly, quantities $\hat{C}_{Y,m'}$,
$\{\hat{\phi}_{Y,m'}^{j}\}^{N_2}_{j=1}$ and
$\{\hat{\lambda}_{Y,m'}^{j}\}^{N_2}_{j=1}$ are defined for the
second sample with $2\leq m'\leq N_2$. To introduce the SN-based
test, we define the recursive estimates
\[
c_{k}^{i,j}=\bigl\langle(\hat{C}_{X, \lfloor kN_1/N \rfloor}-
\hat{C}_{Y,\lfloor
kN_2/N\rfloor})\hat{\phi}_{XY}^{i},\hat{
\phi}_{XY}^{j}\bigr\rangle,\qquad 2\leq k\leq N, 1\leq i,j
\leq K,
\]
which estimate the difference of the covariance operators on the
space spanned by $\{\tilde{\phi}^j\}_{j=1}^K$. Here $K$ is a
user-chosen number, which is held fixed in the asymptotics. Denote
by $\hat{\alpha}_{k}=\operatorname{vech}(\mathbf{C}_k)$ with
$\mathbf{C}_k=(c_{k}^{i,j})_{i,j=1}^{K}$. In the independent and
Gaussian case, Panaretos \textit{et al.} \cite{r19} proposed the following test
(hereafter, the PKM test),
\[
T_{N_1,N_2}=\frac{N_1N_2}{2N}\sum^{K}_{i=1}
\sum^{K}_{j=1}\frac
{(c_{N}^{i,j})^2}{\hat{\varrho}_i\hat{\varrho}_j},\qquad \hat{
\varrho }_j=\frac{1}{N} \Biggl\{\sum
^{N_1}_{i=1}\bigl(\bigl\langle X_i,\hat{
\phi }_{XY}^j\bigr\rangle\bigr)^2+\sum
^{N_2}_{i=1}\bigl(\bigl\langle Y_i,\hat{
\phi}_{XY}^j\bigr\rangle \bigr)^2 \Biggr\},
\]
which converges to $\chi^2_{(K+1)K/2}$ under the null. To take the
dependence into account, we introduce the SN matrix
%
%
\begin{equation}
V_{SN,N}^{(1)}(d)=\frac{1}{N^2}\sum
^{N}_{k=1}k^2(\hat{\alpha
}_{k}-\hat{\alpha}_{N}) (\hat{\alpha}_{k}-\hat{
\alpha}_{N})',
\end{equation}
with $d=(K+1)K/2$. The SN-based test statistic is then defined as,
%
%
\begin{equation}
G_{SN,N}^{(1)}(d)=N\hat{\alpha}_{N}'
\bigl(V_{SN,N}^{(1)}(d)\bigr)^{-1}\hat {
\alpha}_{N}.
\end{equation}
Notice that the PKM test statistic can also be written as a
quadratic form of $\hat{\alpha}_{N}$ but with a different
normalization matrix that is only applicable to the independent and
Gaussian case. The special form of the SN-based test statistic makes
it robust to the dependence within each sample and also the
dependence between the two samples when their sample sizes are
equal. We shall study the asymptotic behavior of $G_{SN,N}^{(1)}(d)$
under the weak dependence assumption in Section~\ref{theory}.

\subsection{Eigenvalues and eigenfunctions}\label{test-eigen-value}

In practice, it is also interesting to infer how far the marginal
distributions of two sequences of stationary functional time series
coincide/differ and quantify the difference. By the
Karhunen--Lo\`{e}ve expansion (Bosq \cite{r4}, page 26), we have
\[
X_i(t)=\sum^{+\infty}_{j=1}\sqrt
{\lambda_X^j}\beta_{X_i,j}\phi
_X^j(t),\qquad Y_i(t)=\sum
^{+\infty}_{j=1}\sqrt{\lambda_Y^j}
\beta _{Y_i,j}\phi_Y^j(t),
\]
where $\beta_{X_i,j}=\int_{\mathcal{I}}X_i(t)\phi_X^j(t)\,\mathrm
{d}t$ and
$\beta_{Y_i,j}=\int_{\mathcal{I}}Y_i(t)\phi_Y^j(t)\,\mathrm{d}t$
are the
principal components (scores), which satisfy that
$E[\beta_{X_i,j}\beta_{X_i,j'}]=\delta_{jj'}$ and
$E[\beta_{Y_i,j}\beta_{Y_i,j'}]=\delta_{jj'}$. The problem is then
translated into testing the equality of the functional principal
components (FPC's) namely the eigenvalues and eigenfunctions. For a
prespecified positive integer $M$, we denote the vector of the first
$M$ eigenvalues by $\lambda_X^{1:M}=(\lambda_X^1,\ldots,\lambda_X^M)$
and $\lambda_Y^{1:M}=(\lambda_Y^1,\ldots,\lambda_Y^M)$. Further
define $\phi_{X}^{1:M}=(\phi_X^1,\ldots,\phi_X^M)$ and
$\phi_{Y}^{1:M}=(\phi_Y^1,\ldots,\phi_Y^M)$ the first $M$
eigenfunctions of the covariance operators $C_X$ and $C_Y$,
respectively. Since the eigenfunctions are determined up to a sign,
we assume that $\langle \phi_X^j,\phi_Y^j\rangle \geq0$ in order for the
comparison to be meaningful. We aim to test the null hypothesis
$H_{2,0}\dvtx \lambda_X^{1:M}=\lambda_Y^{1:M}$ and $H_{3,0}\dvtx
\phi_X^{1:M}=\phi_Y^{1:M}$ versus the alternatives that
$H_{2,a}\dvtx \lambda_X^{1:M}\neq\lambda_Y^{1:M}$ and $H_{3,a}\dvtx
\phi_X^{1:M}\neq\phi_Y^{1:M}$ (in the $L^2$ norm sense). The
problem of comparing the FPC's of two independent and identically
distributed (i.i.d.) functional sequences has been considered in
Benko \textit{et al.} \cite{r2}, where the authors proposed an i.i.d.
bootstrap method
which seems not applicable to the dependent case. The block
bootstrap based method is expected to be valid in the weakly
dependent case but the choice of the block size seems to be a
difficult task in the current setting. To accommodate the dependence
and avoid the bandwidth choice, we adopt the SN idea.

Recall the recursive estimates of the eigenvalues
$\hat{\lambda}^j_{X,m}$ and $\hat{\lambda}^j_{Y,m'}$ which are
calculated based on the subsamples $\{X_i(t)\}_{i=1}^{m}$ and
$\{Y_i(t)\}_{i=1}^{m'}$. Let
$\hat{\theta}_{k}^j=\hat{\lambda}^j_{X,\lfloor
kN_1/N\rfloor}-\hat{\lambda}^j_{Y, \lfloor kN_2/N \rfloor}$ and
$\hat{\theta}_{k}=(\hat{\theta}_{k}^1,\ldots,\hat{\theta}_{k}^M)'$
with $\lfloor N\epsilon\rfloor\leq k\leq N$ for some
$\epsilon\in(0,1]$, which is held fixed in the asymptotics. We
consider the trimmed SN-based test statistic
%
%
\begin{equation}
G_{SN,N}^{(2)}(M)=N^3\hat{\theta}_{N}'
\Biggl\{\sum^{N}_{k=\lfloor
N\epsilon
\rfloor}k^2(\hat{
\theta}_{k}-\hat{\theta}_{N}) (\hat{\theta }_{k}-
\hat{\theta}_{N})' \Biggr\}^{-1}\hat{
\theta}_{N}.
\end{equation}
The trimmed version of the SN-based test statistic is proposed out
of technical consideration when the functional observations lie on
an infinite dimensional space. It can be seen from the proof in the
supplemental material \cite{r25} that the trimming is not required when functional data lie
on a finite-dimensional space; see Remark~0.1
in the
supplemental material \cite{r25}.
%

\begin{remark}
To compare the difference between the eigenvalues, one may also
consider their ratios. Define
$\hat{\zeta}_{k}=(\hat{\lambda}^1_{X,\lfloor
kN_1/N\rfloor}/\hat{\lambda}^1_{Y,\lfloor
kN_2/N\rfloor},\ldots,\hat{\lambda}^M_{X,\lfloor
kN_1/N\rfloor}/\hat{\lambda}^M_{Y,\lfloor kN_2/N\rfloor})'$ for
$k=\lfloor N\epsilon\rfloor,\ldots,N$. An alternative SN-based test
statistic is given by
%
%
\begin{equation}
\tilde{G}_{SN,N}^{(2)}(M)=N(\hat{\zeta}_{N}-
\mathbf{l}_M)' \Biggl\{ \frac{1}{N^2}\sum
^{N}_{k=\lfloor
N\epsilon\rfloor}k^2(\hat{\zeta}_{k}-
\hat{\zeta}_{N}) (\hat{\zeta}_{k}-\hat{\zeta}_{N})'
\Biggr\}^{-1}(\hat{\zeta }_{N}-\mathbf{l}_M),
\end{equation}
where $\mathbf{l}_M$ is a $M$-dimensional vector of all ones. Since
the finite sample improvement by using $\tilde{G}_{SN,N}^{(2)}(M)$
is not apparent, we do not further investigate the properties of
$\tilde{G}_{SN,N}^{(2)}(M)$. 
\end{remark}

We now turn to the problem of testing the equality of the
eigenfunctions. To proceed, we let
%
%
\begin{equation}
\label{basis} \hat{\nu}_j=\bigl(\hat{\phi}_{XY}^{j+1},
\hat{\phi}_{XY}^{j+2},\ldots ,\hat{\phi}_{XY}^p
\bigr)
\end{equation}
be a vector of $p-j$ orthonormal basis functions for $j=1,2,\ldots,M$
with $M\leq p$ and $p$ being a user chosen number. Recall that
$\hat{\phi}^j_{X,m}(t)$ and $\hat{\phi}^j_{Y,m'}(t)$ are the $j$th
eigenfunctions of the empirical covariance operators $\hat{C}_{X,m}$
and $\hat{C}_{Y,m'}$ which are computed based on the first $m$ (and
$m'$) samples. Here we require that
$\langle \hat{\phi}^j_{X,m},\hat{\phi}^j_{X,N_1}\rangle \geq0$ and
$\langle \hat{\phi^j}_{Y,m'},\hat{\phi}^j_{X,N_1}\rangle \geq0$ for $2\leq m\leq
N_1$ and \mbox{$2\leq m'\leq N_2$}. As the eigenfunctions are defined on an
infinite-dimensional space, we project the difference between the
$j$th eigenfunctions onto the space spanned by $\hat{\nu}_j$.
Formally, we define the projection vectors
\[
\hat{\eta}_{k}^{j}=\bigl(\bigl\langle\hat{
\phi}_{X, \lfloor kN_1/N\rfloor}^{j} -\hat{\phi}_{Y,\lfloor
kN_2/N\rfloor}^{j},
\hat{\phi}_{XY}^{j+1}\bigr\rangle,\ldots,\bigl\langle\hat{\phi
}_{X,\lfloor
kN_1/N\rfloor}^{j}-\hat{\phi}_{Y,\lfloor
kN_2/N\rfloor}^{j},\hat{
\phi}_{XY}^{p}\bigr\rangle\bigr),
\]
where $1\leq j\leq M$
and $k=\lfloor N\epsilon\rfloor,\ldots,N$. Further let
$\hat{\eta}_{k}=(\hat{\eta}_{k}^{1},\hat{\eta}_{k}^{2},\ldots
,\hat{\eta}_{k}^{M})'\in\mathbb{R}^{M_0}$
with $M_0=\frac{M(2p-M-1)}{2}$. The trimmed SN-based test statistic
is then defined as
%
%
\begin{equation}
G_{SN,N}^{(3)}(M_0)=N\hat{\eta}_{N}'
\Biggl\{\frac{1}{N^2}\sum^{N}_{k=\lfloor
N\epsilon\rfloor}k^2(
\hat{\eta}_{k}-\hat{\eta}_{N}) (\hat{\eta}_{k}-
\hat{\eta}_{N})' \Biggr\}^{-1}\hat{
\eta}_{N},
\end{equation}
for some $0<\epsilon<1$.
%

\begin{remark}
It is worth noting that $G_{SN,N}^{(3)}(M_0)$ is designed for
testing the equality of the first $M$ eigenfunctions. Suppose we are
interested in testing the hypothesis for a particular eigenfunction,
that is, the null $\phi_X^{j}=\phi_Y^{j}$ versus the alternative
$\phi_X^{j}\neq\phi_Y^{j}$. We can consider the basis functions
\[
\tilde{\nu}_j=\bigl(\hat{\phi}_{XY}^1\ldots,
\hat{\phi }_{XY}^{j-1},\hat{\phi}_{XY}^{j+1},
\ldots,\hat{\phi}_{XY}^p\bigr),
\]
and the projection vector
$\hat{\eta}_{k}^{j}=(\langle\hat{\phi}_{X,\lfloor kN_1/N\rfloor}^{j}
-\hat{\phi}_{Y,\lfloor
kN_2/N\rfloor}^{j},\hat{\phi}_{XY}^1\rangle,\ldots,\langle\hat{\phi
}_{X,\lfloor
kN_1/N\rfloor}^{j} -\hat{\phi}_{Y,\lfloor
kN_2/N\rfloor}^{j},\allowbreak  \hat{\phi}_{XY}^{j-1}\rangle,\langle\hat{\phi}_{X,\lfloor
kN_1/N\rfloor}^{j} -\hat{\phi}_{Y,\lfloor
kN_2/N\rfloor}^{j},\hat{\phi}_{XY}^{j+1}\rangle,\ldots,\langle\hat{\phi
}_{X,\lfloor
kN_1/N\rfloor}^{j}-\hat{\phi}_{Y,\lfloor
kN_2/N\rfloor}^{j},\hat{\phi}_{XY}^{p}\rangle)'$. The SN-based test
statistic can then be constructed in a similar manner. We also note
that when $\phi_X^j\neq\phi_Y^j$ and $\phi_X^i=\phi_Y^i$ for
$i\neq
j$, the choice of $\tilde{\nu}_j$ may result in trivial power
because $\langle \phi_X^j-\phi_Y^j,\tilde{\phi}^i\rangle $ for $i\neq j$ can be
close to 0. In this case, one remedy is to consider alternative
basis functions, for example, (\ref{nu1}) and (\ref{nu2}) as
suggested in
the simulation. 
\end{remark}

%
\begin{remark}\label{nondegenerate}
The choice of the basis functions $\hat{\nu}_j$ is motivated by
the Bahadur representation of the recursive estimates in the
supplemental material \cite{r25}. Under suitable assumptions as given in the next section,
it can be shown that
%
%
\begin{equation}
\label{Bahadur-rep} \bigl\langle \hat{\phi}^a_{X,k},\phi\bigr
\rangle =\bigl\langle \phi^a_X,\phi\bigr\rangle -
\frac{1}{k}\sum^{k}_{i=1} \biggl\{\sum
_{s\neq
a}\frac{\beta_{X_i,s}\beta_{X_i,a}}{\lambda_X^s-\lambda_X^a}\bigl\langle\phi
_X^s,\phi\bigr\rangle \biggr\} +R_{X,k}^a,
\end{equation}
with $R_{X,k}^a$ being the remainder term and $\phi\in
L^2(\mathcal{I})$. The second term on the RHS of (\ref{Bahadur-rep})
plays a key role in determining the limiting distribution of the
SN-based test statistic. When $\phi=\phi_X^j$ with $j\neq a$, the
linear term reduces to
$-\frac{1}{k}\sum^{k}_{i=1}\frac{\beta_{X_i,j}\beta
_{X_i,a}}{\lambda_X^j-\lambda_X^a}$,
which satisfies the functional central limit theorem under suitable
weak dependence assumption. Notice\vadjust{\goodbreak} that the linear term vanishes
when $\phi=\phi_X^a$ and the asymptotic distribution of the
projection vector is degenerate. It is also worth noting that the
linear terms in the Bahadur representations of
$\langle \hat{\phi}^a_{X,k},\phi_X^j\rangle $ and $\langle \hat{\phi}^j_{X,k},\phi_X^a\rangle $
are opposite of each other which suggests that when testing the
eigenfunctions jointly, the basis functions should be chosen in a
proper way so that the asymptotic covariance matrix of the
projection vector, that is, $\hat{\eta}_k$ is nondegenerate. 
\end{remark}

\section{Theoretical results}\label{theory}

To study the asymptotic properties of the proposed statistics, we
adopt the dependence measure proposed in H\"{o}rmann and Kokoszka \cite
{r11}, which is applicable to the temporally dependent functional
process. There are also other weak dependence measures (e.g.,
mixing) or specific processes (e.g., functional linear processes)
suitable for the asymptotic analysis of functional time series (see
Bosq \cite{r4}), we decide to use H\"{o}rmann and Kokoszka's
$L_p$-m-approximating dependence measure for its broad applicability
to linear and nonlinear functional processes as well as its
theoretical convenience and elegance.

%
\begin{definition}\label{def1}
Assume that $\{X_i\}\in L^p_{\mathbb{H}}$ with $p>0$ admits the
following representation
%
%
\begin{equation}
X_i=f(\varepsilon_i, \varepsilon_{i-1},
\ldots),\qquad i=1,2,\ldots,
\end{equation}
where the $\varepsilon_i$'s are i.i.d. elements taking values in a
measurable space $S$ and $f$ is a measurable function $f\dvtx
S^{\infty}\rightarrow\mathbb{H}$. For each $i\in\mathbb{N}$, let
$\{\varepsilon_j^{(i)}\}_{j\in\mathbb{Z}}$ be an independent copy
of $\{\varepsilon_j\}_{j \in\mathbb{Z}}$. The sequence $\{X_i\}$ is
said to be $L^p$-$m$-approximable if
%
%
\begin{equation}
\sum^{\infty}_{m=1}\bigl(E\bigl \Vert
X_m-X_m^{(m)}\bigr \Vert^p
\bigr)^{1/p}<\infty,
\end{equation}
where $X_i^{(m)}=f(\varepsilon_i, \varepsilon_{i-1},\ldots,
\varepsilon_{i-m+1}, \varepsilon_{i-m}^{(i)},
\varepsilon_{i-m-1}^{(i)},\ldots)$.
\end{definition}

Define $B_q(r)$ as a $q$-dimensional vector of independent Brownian
motions. For $\epsilon\in[0,1)$, we let
\[
W_q(\epsilon)=B_q(1)'J_q(
\epsilon)^{-1}B_q(1),\qquad \mbox{where } J_q(
\epsilon)=\int^{1}_{\epsilon}\bigl(B_q(r)-rB_q(1)
\bigr) \bigl(B_q(r)-rB_q(1)\bigr)'\,
\mathrm{d}r.
\]
The critical values of $W_q:=W_q(0)$ have been tabulated by Lobato
\cite{r17}. In general, the quantiles of $W_q(\epsilon)$
can be obtained
via simulation. To derive the asymptotic properties of the proposed
tests, we make the following assumptions.

%
\begin{assumption}\label{m-z1}
Assume $\{X_i(t)\}^{+\infty}_{i=1}\subseteq L^2_{\mathbb{H}}$ and
$\{Y_i(t)\}^{+\infty}_{i=1}\subseteq L^2_{\mathbb{H}}$ are both
$L^4$-$m$-approximable and they are mutually independent.
\end{assumption}

%
\begin{assumption}\label{m-z2}
Assume $\{(X_i(t),Y_i(t))\}^{+\infty}_{i=1}\subseteq
L^4_{\mathbb{H}\times\mathbb{H}}$ is an $L^4$-$m$-approximable
sequence.\vadjust{\goodbreak}
\end{assumption}

%
\begin{assumption}\label{eigen}
Assume $\lambda_X^{1}>\lambda_X^{2}>\cdots>\lambda_X^{m_0+1}$ and
$\lambda_Y^{1}>\lambda_Y^{2}>\cdots>\lambda_Y^{m_0+1}$, for some
positive integer $m_0\geq2$.
\end{assumption}

Note that Assumption~\ref{m-z2} allows dependence between
$\{X_i(t)\}$ and $\{Y_i(t)\}$, which is weaker than Assumption~\ref{m-z1}. To investigate the asymptotic properties of
$G_{SN,N}^{(1)}(d)$ under the local alternatives, we consider the
local alternative $H_{1,a}\dvtx C_X-C_Y=L\bar{C}/\sqrt{N}$ with $\bar{C}$
being a Hilbert--Schmidt operator, where $L$ is a nonzero constant.
Define $\Delta=(\langle\bar{C}\tilde{\phi}^i,\tilde{\phi
}^j\rangle)_{i,j=1}^K\in
\mathbb{R}^{K\times K}$ as the projection of $\bar{C}$ onto the
space spanned by
$\{\tilde{\phi}^1,\tilde{\phi}^2,\ldots,\tilde{\phi}^K\}$ and assume
that $\operatorname{vech}(\Delta)\neq\allowbreak  \mathbf{0}\in\mathbb{R}^{d}$. The
following theorem states the asymptotic behaviors of
$G_{SN,N}^{(1)}(d)$ under the null and the local alternatives.

%
\begin{theorem}\label{SN-cov}
Suppose Assumptions~\ref{m-z1}, \ref{eigen} hold with $m_0\geq K$.
Further assume that the asymptotic covariance matrices
$\Lambda^*_{d}(\Lambda_d^*)'$ given in Lemma~0.3
is
positive definite. Then under $H_{1,0}$,
$G_{SN,N}^{(1)}(d)\rightarrow^d W_d$ and under $H_{1,a}$,
$\lim_{|L|\rightarrow+\infty}\lim_{N\rightarrow+\infty
}G_{SN,N}^{(1)}(d)=+\infty$.
Furthermore, if $\gamma_1=\gamma_2$, then the conclusion also holds
with Assumption~\ref{m-z1} replaced by Assumption~\ref{m-z2}.
\end{theorem}

It is seen from Theorem~\ref{SN-cov} that $G_{SN,N}^{(1)}(d)$ has
pivotal limiting distributions under the null and they are
consistent under the local alternatives as $L\rightarrow+\infty$. It
is worth noting that in our asymptotic framework, $d$ (or $K$) is
assumed to be fixed as $n\rightarrow\infty$. Since $K$ is usually
chosen to make the first $K$ principle components explain a certain
percentage of variation (say $85\%$), the magnitude of $K$
critically depends on the prespecified threshold and the decay rate
of the eigenvalues. In some cases, $d=(K+1)K/2$ can be quite large
relative to sample size so it may be more meaningful to use the
asymptotic results established under the framework that
$d\rightarrow\infty$ but $d/n\rightarrow0$ as
$n\rightarrow\infty$. This motivates the question that whether the
following convergence result
\[
\sup_{x\in\R}\bigl |P\bigl(G_{SN,N}^{(1)}(d)\le x
\bigr)-P(W_d\le x)\bigr |\rightarrow0 \qquad\mbox{as } n\rightarrow\infty
\]
holds
when $d$ diverges to $\infty$ but at a slower rate than $n$. This
would be an interesting future research topic but is beyond the
scope of this paper.

To study the asymptotics of $G_{SN,N}^{(2)}(M)$ and
$G_{SN,N}^{(3)}(M_0)$, we introduce some notation. Let
$\omega_{X_i}^{jk}=\beta_{X_i,j}\beta_{X_i,k}$ and
$r_{X}^{jk,j'k'}(h)=E[(\omega_{X_i}^{jk}-\delta_{jk}\lambda
_j)(\omega_{X_{i+h}}^{j'k'}-\delta_{j'k'}\lambda_{j'})]$
be the cross-covariance function between $\omega_{X_i}^{jk}$ and
$\omega_{X_i}^{j'k'}$ at lag $h$. Set $r^{jk}_X(h):=r^{jk,jk}_X(h)$.
Define
$v_{X_i}^{jk}=\omega_{X_i}^{jk}-E[\omega_{X_i}^{jk}]=\omega
_{X_i}^{jk}-\delta_{jk}\lambda_j$.
Analogous\vspace*{1.5pt} quantities $r_{Y}^{jk,j'k'}(h)$ and $v_{Y_i}^{jk}$ can be
defined for the second sample. We make the following assumption to
facilitate our derivation.
%

\begin{assumption}\label{long-run}
Suppose that
%
%
\begin{equation}
\sum_{j,k}\sum_{j',k'}
\Biggl(\sum^{+\infty}_{h=-\infty
}\bigl |r^{jk,j'k'}_X(h)\bigr |
\Biggr)^2<+\infty,\qquad \sum_{j,k}\sum
^{+\infty}_{h=-\infty}\bigl |r^{jk}_X(h)\bigr |<+
\infty
\end{equation}
and
%
%
\begin{equation}
\sum_{j,k}\sum_{j',k'}
\sum_{i_1,i_2,i_3\in\mathbb{Z}}\bigl |\operatorname{cum}\bigl(v_{X_0}^{jk},v_{X_{i_1}}^{jk},v_{X_{i_2}}^{j'k'},v_{X_{i_3}}^{j'k'}
\bigr)\bigr |<\infty .
\end{equation}
The summability conditions also hold for the second sample
$\{Y_i(t)\}$.
\end{assumption}

Assumption~\ref{long-run} is parallel to the summability condition
considered in Benko \textit{et al.} \cite{r2} (see Assumption~1
therein) for i.i.d.
functional data. It is not hard to verify the above assumption for
Gaussian linear functional process (see, e.g., Bosq
\cite{r4}), as
demonstrated in the following proposition.

%
\begin{proposition}\label{glp}
Consider the linear process
$X_i(t)=\sum^{\infty}_{j=0}b_j\varepsilon_{i-j}(t)$, where
$\varepsilon_j(t)=   \sum^{\infty}_{i=1}\sqrt{\lambda_i}z_{i,j}\phi_i(t)$
with $\{z_{i,j}\}$ being a sequence of independent standard normal
random variables across both index $i$ and $j$. Let
$\pi(h)=\sum_{i}b_ib_{i+h}$. Assume that
$\sum^{\infty}_{j=1}\lambda_j<\infty$ and $\sum_h|\pi
(h)|<\infty$.
Then Assumption~\ref{long-run} holds for $\{X_i(t)\}$.
\end{proposition}

%
\begin{theorem}\label{first-value}
Suppose Assumptions~\ref{m-z1}, \ref{eigen}, \ref{long-run} hold
with $m_0 \geq M$ and the asymptotic covariance matrix
$\tilde{\Lambda}_M\tilde{\Lambda}_M'$ given in Lemma~0.5
is positive definite. Then under $H_{2,0}$, we have
$G_{SN,N}^{(2)}(M)\rightarrow^d W_M(\epsilon)$. Under the local
alternative
$H_{2,a}\dvtx \lambda_X^{1:M}-\lambda_Y^{1:M}=\frac{L}{\sqrt{N}}\bar
{\lambda}$
with $\bar{\lambda}\neq \mathbf{0}\in\mathbb{R}^M$, we have
$\lim_{|L|\rightarrow\infty}\lim_{N\rightarrow
+\infty}G_{SN,N}^{(2)}(M)=+\infty$.
\end{theorem}

In order to study the asymptotic properties of $G_{SN,N}^{(3)}(M_0)$
under the null and local alternative, we further make the following
assumption.

%
\begin{assumption}\label{assefc}
Let $\beta_{X_i,j}^{(m)}=\int X_i^{(m)}(t)\phi_{X}^j(t)\,\mathrm
{d}t$, where
$X_i^{(m)}$ is the $m$-dependent approximation of $X_i(t)$ (see
Definition~\ref{def1}). Suppose one of the following conditions
holds:
%
%
\begin{equation}
\label{addsum} \sum^{\infty}_{m=1}\sum
_{
j=1}^{\infty} \bigl\{E \bigl(\beta_{X_1,j}-
\beta_{X_1,j}^{(m)} \bigr)^4 \bigr\}^{1/4}<
\infty,\qquad \sum_{j=1}^{\infty} \bigl(E
\beta_{X_1,j}^4 \bigr)^{1/4}<\infty,
\end{equation}
or
%
%
\begin{equation}
\label{addsum2} \sum^{+\infty}_{s=1}\bigl |\bigl\langle
\phi_X^s,\tilde{\phi}^j\bigr\rangle \bigr |<+
\infty,\qquad 2\leq j\leq p.
\end{equation}
The same condition holds for the second sample $\{Y_i(t)\}$.
\end{assumption}

%
\begin{theorem}\label{eigenfunction}
Suppose Assumptions~\ref{m-z1}, \ref{eigen}, \ref{long-run} and~\ref{assefc} hold with $m_0 \geq M$ and the asymptotic covariance
matrix $\bar{\Lambda}_{M_0}\bar{\Lambda}_{M_0}'$ given in Lemma~0.7
is positive definite. Then under $H_{3,0}$, we have
$G_{SN,N}^{(3)}(M_0)\rightarrow^d W_{M_0}(\epsilon)$.
\end{theorem}

%
\begin{proposition}\label{alternative}
Define $\tilde{\Delta}$ by replacing $\hat{\phi}_{X,N_1}^{j}$,
$\hat{\phi}_{Y,N_2}^{j}$ and $\hat{\phi}_{XY}^j$ with $\phi_X^{j}$,
$\phi_Y^{j}$ and $\tilde{\phi}^j$ in the definition of
$\hat{\eta}_{N}$. Consider the local alternative $H_{3,a}\dvtx
\tilde{\Delta}=L\boldsymbol{\bar{\psi}}/\sqrt{N}$ with
$\boldsymbol{\bar{\psi}}\neq\mathbf{0}\in\mathbb{R}^{M_0}$. Suppose
Assumptions~\ref{m-z1}, \ref{eigen}, \ref{long-run} and~\ref{assefc}
hold with $m_0\geq M$ and the asymptotic covariance matrix
$\bar{\Lambda}_{M_0}\bar{\Lambda}_{M_0}'$ given in Lemma~0.7
is positive definite. Then we have
\[
\lim_{|L|\rightarrow\infty}\lim_{N\rightarrow+\infty
}G_{SN,N}^{(3)}(M_0)=+\infty
\]
under $H_{3,a}$.
\end{proposition}

It is worth noting that the conclusions in Theorem~\ref{first-value}, Theorem~\ref{eigenfunction} and Proposition~\ref{alternative} also hold with Assumption~\ref{m-z1} replaced by
Assumption~\ref{m-z2} and $\gamma_1=\gamma_2$. Finally, we point out
that condition (\ref{addsum}) can be verified for Gaussian linear
functional process as shown in the following proposition.
%

\begin{proposition}\label{glp2}
Consider the Gaussian linear process in Proposition~\ref{glp}.
Assume that $\sum^{\infty}_{j=1}\sqrt{\lambda_j}<\infty$ and
$\sum^{\infty}_{m=1}(\sum^{\infty}_{j=m}b_j^2)^{1/2}<\infty$. Then
Assumption~\ref{long-run} and condition (\ref{addsum}) are satisfied
for $\{X_i(t)\}$.
\end{proposition}

\section{Numerical studies}\label{simulation}

We conduct a number of simulation experiments to assess the
performance of the proposed SN-based tests in comparison with the
alternative methods in the literature. We generate functional data
on a grid of $10^3$ equispaced points in $[0,1]$, and then convert
discrete observations into functional objects by using B-splines
with 20 basis functions. We also tried 40 and 100 basis functions
and found that the number of basis functions does not affect our
results much. Throughout the simulations, we set the number of Monte
Carlo replications to be 1000 except for the i.i.d. bootstrap method in
Benko \textit{et al.} \cite{r2}, where the number of replications is
only 250
because of high computational cost.

\begin{table}[b]
\tablewidth=\textwidth
\tabcolsep=0pt
\caption{Empirical sizes and size-adjusted powers of (i) the
SN-based test, (ii) the PKM test and (iii) the CLRV
test for testing the equality of the covariance operators. The
nominal level is 5\%}\label{tb-cov}
\begin{tabular*}{\textwidth}{@{\extracolsep{\fill}}lllld{2.1}d{2.1}d{2.1}d{3.1}d{3.1}d{2.1}d{3.1}@{}}
\toprule
&&&&\multicolumn{7}{l@{}}{$K$}
\\[-5pt]
&&&&\multicolumn{7}{l@{}}{\hrulefill}
\\
& Parameter & $N_1=N_2$ & & \multicolumn{1}{l}{1} & \multicolumn{1}{l}{2} & \multicolumn{1}{l}{3} & \multicolumn{1}{l}{4} & \multicolumn{1}{l}{5} & \multicolumn{1}{l}{$K^*_1$} & \multicolumn{1}{l@{}}{$K_2^*$}
\\
\midrule
A & $\mathbf{v}_X=(12,7,0.5,9,5,0.3)$ & 100 & (i) & 4.3 & 5.7 & 6.8 & 8.7 & 14.3  & 8.7 & 10.7\\
&&& (ii) & 14.5 & 20.9 & 22.9 & 32.2 & 39.5 & 32.0 & 34.0 \\
&&& (iii) & 9.1 & 12.9 & 20.9 & 39.8 & 67.9 & 38.8 & 48.6\\
& $\mathbf{v}_Y=(12,7,0.5,9,5,0.3)$ & 200& (i) & 4.7 &  5.7 & 4.6 &  7.0 & 8.0 & 7.0 & 7.0\\
&&& (ii) & 12.8 & 20.6 & 26.7 & 34.7 & 42.6  & 34.8 & 37.3 \\
&&& (iii) & 6.9 & 9.6 & 14.5 & 25.2 & 41.5 & 25.1 & 28.9\\
B & $\mathbf{v}_X=(14,7,0.5,6,5,0.3)$ & 100 &(i)& 19.1 & 23.6 & 17.7 & 14.2 & 12.7 & 14.1 & 13.0\\
&&& (ii)& 27.6 & 37.7 & 31.6 & 22.9 & 21.2  & 23.1 & 22.7\\
&&& (iii) &27.0 & 33.8 & 23.0 & 20.5 & 14.0  & 20.5 & 15.8\\
& $\mathbf{v}_Y=(8,7,0.5,6,5,0.3)$ & 200 & (i)& 31.2 & 37.7 & 30.4 & 21.9 & 21.9 & 21.9 & 22.1\\
&&& (ii)& 39.1 & 61.6 & 51.7 & 44.2 & 41.2 & 44.2 & 41.5\\
&&& (iii) & 37.6 &  57.0 & 44.7 & 30.1 & 24.3 & 30.1 & 24.9\\
C & $\mathbf{v}_X=(12,7,0.5,9,3,0.3)$ & 100 &(i)& 5.5 & 10.9 & 30.8 & 62.4 & 64.7  & 62.3 & 63.7\\
&&& (ii) & 4.7 & 16.1 & 57.3 & 94.6 & 98.7 & 94.4 & 97.1\\
&&& (iii) & 5.5 & 13.4 &42.3 & 79.4 & 70.8 & 79.2 & 74.3 \\
& $\mathbf{v}_Y=(12,7,0.5,3,9,0.3)$ & 200 &(i)& 5.3 & 10.0 & 45.7 & 90.3 & 94.3 & 90.4 & 92.5\\
&&& (ii) &6.4 & 13.0 & 67.8 & 99.9 & 100.0 &  99.9 & 100.0\\
&&& (iii) &6.1 & 12.5 & 60.5 & 99.9 & 99.8 & 99.9 & 99.8\\
D & $\mathbf{v}_X=(12,7,0.5,9,5,0.3)$ & 100 &(i)& 6.1 & 8.3 & 28.3 & 80.1 & 82.2 & 75.6 & 80.6\\
&&& (ii)& 5.5  & 14.6 &47.2 & 100.0 & 100.0 & 94.7 & 100.0\\
&&& (iii) &6.9 & 12.3 & 37.2 & 95.7 & 88.6 & 90.6 & 90.7\\
& $\mathbf{v}_Y=(12,7,0.5,0,5,0.3)$ & 200 &(i)& 5.7 & 8.9 & 39.7 & 96.3 & 98.4 & 95.5 & 98.3 \\
&&& (ii)& 6.4  & 14.5 & 53.6 & 100.0 & 100.0 & 99.4 & 100.0\\
&&& (iii)& 6.0 & 12.9 & 47.7 & 100.0 & 100.0 & 99.3 & 100.0\\
\midrule
\end{tabular*}
\legend{Note: Under the alternatives, we simulate the size-adjusted critical
values by assuming that both $\{X_i\}$ and $\{Y_i\}$ are generated
from (\ref{AR}) with $\rho=0.5$, $\mu=1$ and
$\mathbf{v}=\mathbf{v}_X$.}
\end{table}

\subsection{Comparison of covariance operators}

To investigate the finite sample properties of $G_{SN,N}^{(1)}(d)$
for dependent functional data, we modify the simulation setting
considered in Panaretos \textit{et al.} \cite{r19}. Formally, we consider the
model,
%
%
\begin{equation}
\label{AR} \sum^{3}_{j=1} \bigl\{
\xi_{j,1}^i\sqrt{2}\sin(2\uppi jt)+\xi_{j,2}^i
\sqrt{2}\cos(2\uppi jt) \bigr\},\qquad i=1,2,\ldots, t\in[0,1],
\end{equation}
where the coefficients
$\xi_i=(\xi_{1,1}^i,\xi_{2,1}^i,\xi_{3,1}^i,\xi_{1,2}^i,\xi
_{2,2}^i,\xi_{3,2}^i)'$
are generated from a VAR process,
%
%
\begin{equation}
\label{var} \xi_i=\rho\xi_{i-1}+\sqrt{1-
\rho^2}e_i,
\end{equation}
with $e_i\in\mathbb{R}^6$ being a sequence of i.i.d. normal random
variables with mean zero and covariance matrix
$\Sigma_e=\frac{1}{1+\mu^2}\operatorname{diag}(\mathbf{v})+\frac{\mu
^2}{1+\mu^2}\mathbf{1}_6\mathbf{1}_6'$.
We generate two independent functional time series $\{X_i(t)\}$ and
$\{Y_i(t)\}$ from (\ref{AR}) with $\rho=0.5$ and $\mu=1$. We compare
the SN-based test with the PKM test which is designed for
independent Gaussian process, and the traditional test which is
constructed based on a consistent LRV estimator (denoted by CLRV),
that is,
$G_{CL,N}(d)=N\hat{\alpha}_N\hat{\Sigma}_{\alpha}^{-1}\hat{\alpha}_N$,
where $\hat{\Sigma}_{\alpha}$ is a lag window LRV estimator with
Bartlett kernel and data dependent bandwidth (see Andrews (1991)).
We report the simulation results for $N_1=N_2=100,200$,
$K=1,2,3,4,5$ ($d=1,3,6,10,15$) and various values of $\mathbf{v}$
in Table~\ref{tb-cov}. Results in scenario A show that the size
distortion of all the three tests increases as $K$ gets larger. The
SN-based test has the best size compared to the other two tests. The
PKM test is severely oversized due to the fact that it does not take
the dependence into account. It is seen from the table that the CLRV
test also has severe size distortion especially for large $K$, which
is presumably due to the poor estimation of the LRV matrix of
$\hat{\alpha}_N$ when the dimension is high. Under the alternatives,
we report the size-adjusted power which is computed using finite
sample critical values based on the simulation under the null model
where we assume that both $\{X_i(t)\}$ and $\{Y_i(t)\}$ are
generated from (\ref{AR}) with $\rho=0.5$, $\mu=1$ and
$\mathbf{v}=\mathbf{v}_X$. From scenarios B--D in Table~\ref{tb-cov},
we observe that the PKM is most powerful which is largely due to its
severe upward size distortion. The SN-based test is less powerful
compared to the other two tests but the power loss is generally
moderate in most cases. Furthermore, we present the results when
choosing $K$ by
%
%
\begin{equation}
\label{auto} K^*_j=\operatorname{arg\,min} \biggl\{1\leq J\leq20
\dvtx \frac{\sum^{J}_{i=1}\hat{\lambda}^{i}_{XY}}{\sum^{20}_{i=1}\hat
{\lambda}^i_{XY}}>\alpha_j^* \biggr\},\qquad j=1,2,
\end{equation}
where $\alpha_1^*=85\%$ and $\alpha_2^*=95\%$. An alternative way of
choosing $K$ is to consider the penalized fit criteria (see
Panaretos \textit{et al.} \cite{r19} for the details). We notice that the
performance of all the three tests based on automatic choice $K^*_j$
is fairly close to the performance when $K=4$ or $5$ in most cases.
To sum up, the SN-based test provides the best size under the null
and has reasonable power under different alternatives considered
here, which is consistent with the ``better size but less power''
phenomenon seen in the univariate setup (Lobato \cite
{r17} and Shao \cite{r24}).

\subsection{Comparison of eigenvalues and eigenfunctions}

In this subsection, we study the finite sample performance of the
SN-based test for testing the equality of the eigenvalues and
eigenfunctions. We consider the data generating process,
%
%
\begin{equation}
\label{AR2} \sum^{2}_{j=1} \bigl\{
\xi_{j,1}^i\sqrt{2}\sin(2\uppi jt+\delta_j)+
\xi_{j,2}^i\sqrt{2}\cos(2\uppi jt+\delta_j)
\bigr\},\qquad i=1,2,\ldots, t\in[0,1],
\end{equation}
where $\xi_i^*=(\xi_{1,1}^i,\xi_{2,1}^i,\xi_{1,2}^i,\xi_{2,2}^i)'$
is a 4-variate VAR process (\ref{var}) with $e_i\in\mathbb{R}^4$
being a sequence of i.i.d. normal random variables with mean zero and
covariance matrix
$\Sigma_e=\frac{1}{1+\mu^2}\operatorname{diag}(\mathbf{v})+\frac{\mu
^2}{1+\mu^2}\mathbf{1}_4\mathbf{1}_4'$.
We set $\rho=0.5$ and $\mu=0$. Under $H_{2,0}$ (or $H_{2,a}$),
$\{X_i(t)\}$ and $\{Y_i(t)\}$ are generated independently from
(\ref{AR2}) with $\delta_1=\delta_2=0$ and
$\mathbf{v}_X=\mathbf{v}_Y$ (or $\mathbf{v}_X\neq\mathbf{v}_Y$).
Notice that the eigenvalues of $\{X_i(t)\}$ and $\{Y_i(t)\}$ are
given respectively, by $\mathbf{v}_X$ and $\mathbf{v}_Y$ when
$\delta_1=\delta_2=0$. Under $H_{3,0}$ and $H_{3,a}$, we generate
$\{X_i(t)\}$ and $\{Y_i(t)\}$ independently from (\ref{AR2}) with
$\mathbf{v}_X=\mathbf{v}_Y$, $\delta_{X,1}-\delta_{Y,1}=\delta$, and
$\delta_{X,2}=\delta_{Y,2}=0$, where $\delta=0$ under the null and
$\delta\neq0$ under the alternatives. We aim to test the equality
of the first four eigenvalues and eigenfunctions separately and
jointly. Because functional data are finite dimensional, we
implement the untrimmed version of the SN-based tests, that is,
$\epsilon=0$. To further assess the performance of the SN-based
test, we compare our method with the subsampling approach with
several choices of subsampling widths and the i.i.d. bootstrap method
in Benko \textit{et al.} \cite{r2}. Suppose $N_1=N_2=N_0$. Let $l$ be the
subsampling width and
$\lambda_{\mathrm{sub},i}^j=\lambda_{\mathrm{sub},X,i}^j-\lambda_{\mathrm{sub},Y,i}^j$,
$i=1,2,\ldots, s_{N_0}(l)=\lceil N_0/l \rceil$, where
$\lambda_{\mathrm{sub},X,i}^j$ and $\lambda_{\mathrm{sub},Y,i}^j$ are estimates of the
$j$th eigenvalues based on the $i$th nonoverlapping subsamples
$\{X_k(t)\}^{il}_{k=(i-1)l+1}$ and $\{Y_k(t)\}^{il}_{k=(i-1)l+1}$, respectively. The subsampling variance estimate is given by
$\sigma^2_{\mathrm{sub},j}=\frac{l}{s_{N_0}(l)}\sum^{s_{N_0}(l)}_{i=1}
(\hat{\lambda}
_{\mathrm{sub},i}^j-\frac{1}{s_{N_0}(l)}\sum^{s_{N_0}(l)}_{i=1}\hat{\lambda
}_{\mathrm{sub},i}^j )^2$,
and the test statistic based on the subsampling variance estimate
for testing the equality of the $j$th eigenvalue is defined as
$G_{\mathrm{sub},N}=N_0(\hat{\lambda}_{X,N_0}^j-\hat{\lambda
}_{Y,N_0}^j)^2/\sigma^2_{\mathrm{sub},j}$.
Since the data-dependent rule for choosing the subsampling width is
not available in the current setting, we tried $l=8,12,16$ for
$N_0=48,96$. For testing the equality of eigenvalues jointly and
equality of the eigenfunctions, we shall consider a multivariate
version of the subsampling-based test statistic which can be defined
in a similar fashion. Table~\ref{tb-eigenvalue} summarizes some
selective simulation results for testing the eigenvalues with
various values of $\mathbf{v}$. From scenario A, we see that
performance of the SN-based test under the null is satisfactory
while the size distortion of the subsampling-based method is quite
severe and is sensitive to the choice of block size $l$. It is also
not surprising to see that the i.i.d. bootstrap method has obvious size
distortion as it does not take the dependence into account. Under
the alternatives (scenarios B--D), we report the size-adjusted power
by using the simulated critical values as described in previous
subsection. When the sample size is 48, the SN-based method delivers
the highest power among the tests and it tends to have some moderate
power loss when the sample size increases to 96. On the other hand,
the subsampling method is sensitive to the choice of subsampling
width and its power tends to decrease when a larger subsampling
width is chosen.

\begin{table}
\tablewidth=\textwidth
\tabcolsep=0pt
\caption{Empirical sizes and size-adjusted powers of (i) the
SN-based test, the subsampling-based test with (ii) $l=8$,
(iii) $l=12$ and (iv) $l=16$, and (v) Benko et
al.'s i.i.d. bootstrap based method for testing the equality of the
first two eigenvalues separately (the columns with $M=1,2$) and
jointly (the column with $M=(1,2)$), and the equality of the first
four eigenvalues jointly (the column with $M=(1,2,3,4)$). The
nominal level is 5\% and the number of replications for i.i.d.
bootstrap method is 250}\label{tb-eigenvalue}
\begin{tabular*}{\textwidth}{@{\extracolsep{\fill}}lllld{2.1}d{2.1}d{2.3}d{2.6}@{}}
\toprule
&&&&\multicolumn{4}{l@{}}{$M$}
\\[-5pt]
&&&&\multicolumn{4}{l@{}}{\hrulefill}
\\
& Parameter & $N_1=N_2$ & & \multicolumn{1}{l}{1} & \multicolumn{1}{l}{2} & \multicolumn{1}{l}{(1, 2)} & \multicolumn{1}{l@{}}{(1, 2, 3, 4)}
\\
\midrule
A & $\mathbf{v}_X=(10,0.5,5,0.3)$ & 48 & (i) & 5.4 & 5.1 & 4.6 & 3.8 \\
&&& (ii) &24.2 &38.5 &52.4 & 90.8 \\
&&& (iii) &21.9 &28.8 &51.3 & 68.8 \\
&&& (iv) &21.8 &28.1 &57.9 & 44.7 \\
&&& (v) &11.2 &9.2 &11.6  & 11.6\\
& $\mathbf{v}_Y=(10,0.5,5,0.3)$ & 96 & (i) & 5.2 & 5.6 & 4.8 & 5.1  \\
&&& (ii) & 19.0 & 40.4 & 46.4 & 84.8 \\
&&& (iii) & 16.3 & 29.6 & 38.2 & 77.0   \\
&&& (iv) &16.0 & 25.3 & 36.5 & 78.5  \\
&&& (v) &14.4 & 8.4 & 15.2 & 15.2  \\
B & $\mathbf{v}_X=(20,0.5,5,0.3)$ & 48 &(i) & 25.1 & 4.3 & 21.8 & 15.5 \\
&&& (ii) &24.2 & 5.4 & 13.3 & 7.1 \\
&&& (iii) &19.8 & 6.8 & 8.8 & 8.8  \\
&&& (iv) &14.1 & 6.8 & 8.0 & 9.1  \\
& $\mathbf{v}_Y=(10,0.5,5,0.3)$ & 96 & (i)& 48.4 & 4.8 & 35.6 & 25.0  \\
&&& (ii)& 58.4 & 6.9 & 29.4 & 11.4  \\
&&& (iii) & 50.9 &  6.1 & 29.7 & 13.3  \\
&&& (iv) &53.8 &  6.0 & 29.3 & 11.4  \\
C & $\mathbf{v}_X=(10,0.5,5,0.3)$ & 48 &(i)& 6.2 & 70.6 & 58.9 & 44.1  \\
&&& (ii)& 5.5 & 68.1 & 54.6 & 13.9  \\
&&& (iii)& 4.8 & 49.3 & 23.0 & 16.9   \\
&&& (iv)& 6.1 & 34.2 & 15.4 & 21.2   \\
& $\mathbf{v}_Y=(10,0.5,1,0.3)$ & 96 &(i)& 4.7 & 91.4 & 84.6 & 77.6  \\
&&& (ii)& 4.7 & 98.7 & 96.3 & 69.7  \\
&&& (iii)& 5.5 & 97.9 & 92.5 & 51.6   \\
&&& (iv)& 5.4 & 96.5 & 83.0 & 29.4   \\
D & $\mathbf{v}_X=(20,0.5,5,0.3)$ & 48 &(i)& 27.0 & 70.1 & 68.4  & 55.3 \\
&&& (ii)& 25.8 & 65.7 & 51.1 & 14.0  \\
&&& (iii)& 20.9 & 58.4 & 23.9 & 19.9   \\
&&& (iv)& 14.9 & 40.1 & 11.8 & 17.2   \\
& $\mathbf{v}_Y=(10,0.5,1,0.3)$ & 96 &(i)& 55.3 & 87.9 & 88.4  & 83.3 \\
&&& (ii)& 54.5 & 98.3 & 96.9 & 62.3  \\
&&& (iii)& 48.8 & 97.6 & 95.2 & 53.2   \\
&&& (iv)& 50.1 & 95.2 & 88.0 & 27.7   \\
\midrule
\end{tabular*}
\legend{Note: Under the alternatives, we simulate the size-adjusted critical
values by assuming that both $\{X_i\}$ and $\{Y_i\}$ are generated
from (\ref{AR}) with $\rho=0.5$, $\mu=0$ and
$\mathbf{v}=\mathbf{v}_X$.}
\end{table}

\begin{figure}

\includegraphics{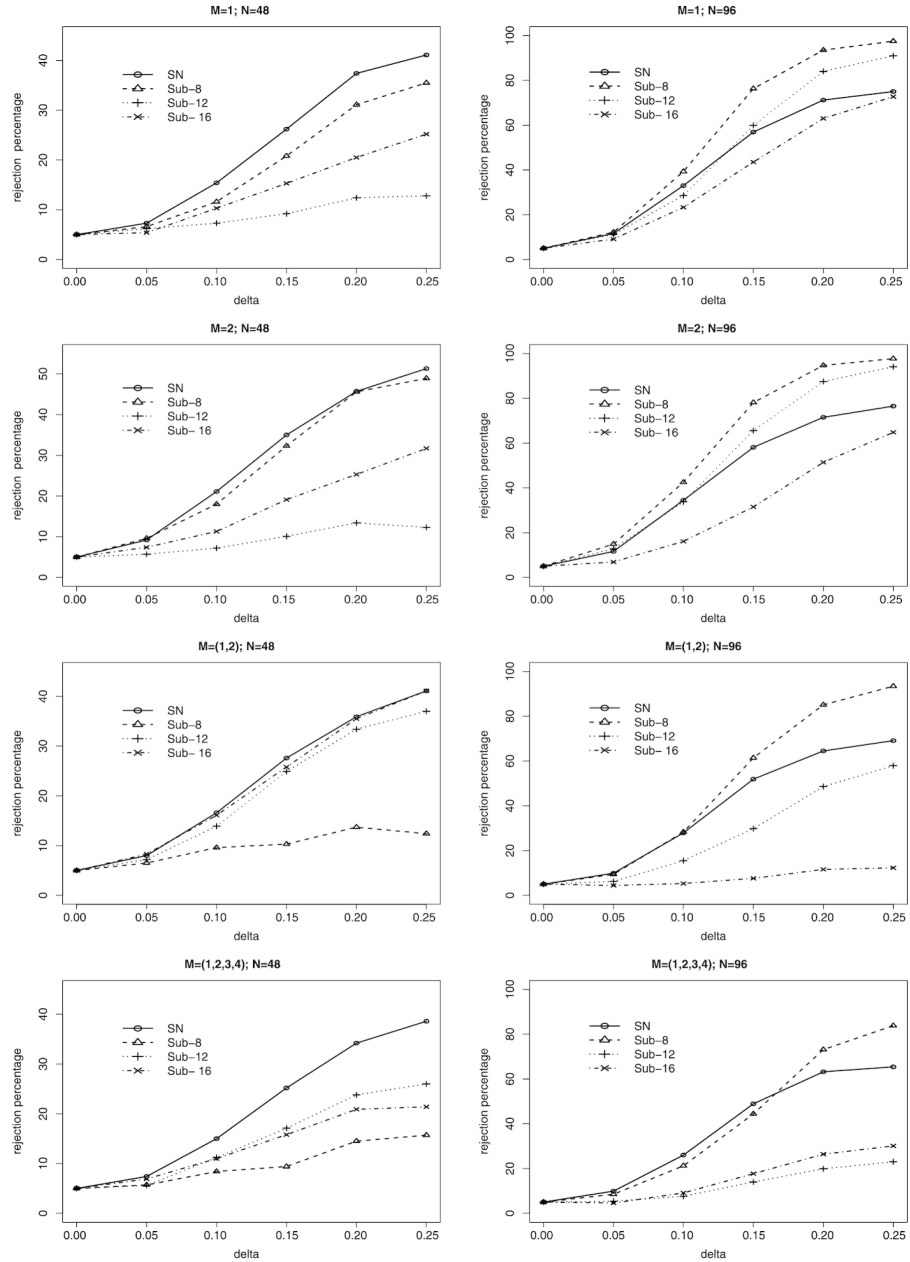}

\caption{Size-adjusted powers of the SN-based test and the
subsampling-based tests for testing the equality of the first two
eigenfunctions separately and jointly, and the equality of the first
four eigenfunctions jointly.}\label{fig:eigenfc-power}
\end{figure}

To test the equality of the first four eigenfunctions, we implement
the SN-based test and the subsampling-based test with the basis
functions,
%
%
\begin{eqnarray}
\label{nu1} &&\hat{\nu}_j^*=(\hat{\phi}_{XY}^1+
\hat{\phi}_{XY}^j,\ldots,\hat {\phi}_{XY}^{j-1}+
\hat{\phi}_{XY}^j,
\nonumber
\\[-8pt]
\\[-8pt]
&&\phantom{\hat{\nu}_j^*=\bigl(}\hat{\phi}_{XY}^{j+1}+
\hat {\phi}_{XY}^j,\ldots,\hat{\phi}_{XY}^p+
\hat{\phi}_{XY}^j\bigr),\qquad 1\leq j\leq4, p=4,
\nonumber
\end{eqnarray}
for testing individual eigenfunction and
%
%
\begin{equation}
\label{nu2} \hat{\nu}_j^{**}=\bigl(\hat{
\phi}_{XY}^{j+1}+\hat{\phi}_{XY}^{j},\hat
{\phi}_{XY}^{j+2}+\hat{\phi}_{XY}^{j},
\ldots,\hat{\phi }_{XY}^p+\hat{\phi}_{XY}^{j}
\bigr),\qquad 1\leq j\leq j^*, p=4,
\end{equation}
with $j^*=2,3,4$, for testing the first $j^*$ eigenfunctions jointly
(correspondingly $M_0=3,5,6$). The tests with the above basis
functions tend to provide similar sizes but higher powers as
compared to the tests with the basis functions $\hat{\nu}_i$ in our
simulation study. The basis\vspace*{-2pt} functions $\hat{\nu}_j^*$ is constructed
by adding the same estimated eigenfunction $\hat{\phi}_{XY}^{j}$ to
each component of $\hat{\nu}_j$, and the associated SN-based test is
expected to be asymptotically valid in view of the Bahadur
representation (\ref{Bahadur-rep}). Selective simulation results are
summarized in Table~\ref{tb-fcn} and Figure~\ref{fig:eigenfc-power}
which present the sizes of the SN-based test, the subsampling-based
test and the i.i.d. bootstrap method, and the size adjusted powers of
the former two respectively. It is seen from Table~\ref{tb-fcn} that
the sizes of the SN-based test are accurate while the
subsampling-based test is apparently size-distorted. It is somewhat
surprising to see that the i.i.d. bootstrap provides better sizes
compared to the subsampling-based approach which is designed for
dependent data. Figure~\ref{fig:eigenfc-power} plots the
(size-adjusted) power functions of the SN-based test and the
subsampling-based test which are monotonically increasing on
$\delta$. When $N_1=N_2=48$, the SN-based test delivers the highest
power in most cases. The subsampling-based test with a small
subsampling width becomes most powerful when sample size increases
to 96. Overall, the SN-based test is preferable as it provides quite
accurate size under the null and has respectable power under the
alternatives.

\begin{table}
\tablewidth=\textwidth
\tabcolsep=0pt
\caption{Empirical sizes of (i) the SN-based test, the
subsampling-based test with (ii) $l=8$, (iii) $l=12$
and (iv) $l=16$, and (v) Benko et al.'s i.i.d. bootstrap
based method for testing the equality of the first two
eigenfunctions separately (the columns with $M=1,2$) and jointly
(the column with $M=(1,2)$), and the equality of the first four
eigenfunctions jointly (the column with $M=(1,2,3,4)$). The nominal
level is 5\% and the number of replications for i.i.d. bootstrap is
250}\label{tb-fcn}
\begin{tabular*}{\textwidth}{@{\extracolsep{\fill}}lllld{2.1}d{2.1}d{2.3}d{2.6}@{}}
\toprule
&&&&\multicolumn{4}{l@{}}{$M$}
\\[-5pt]
&&&&\multicolumn{4}{l@{}}{\hrulefill}
\\
& Parameter & $N_1=N_2$ & & \multicolumn{1}{l}{1} & \multicolumn{1}{l}{2} & \multicolumn{1}{l}{(1, 2)} & \multicolumn{1}{l@{}}{(1, 2, 3, 4)}
\\
\midrule
A & $\mathbf{v}_X=(10,0.5,5,0.3)$ & 48 & (i) &6.4 &3.2 &4.9 & 4.8  \\
&&& (ii) &39.5 &36.0 &78.7 & 67.0 \\
&&& (iii) &62.9 &62.3 &24.8 & 26.4 \\
&&& (iv) &32.6 &25.9 &9.6  & 7.7\\
&&& (v) &2.4 &11.2  &2.4 &  2.4\\
& $\mathbf{v}_Y=(8,0.5,4,0.3)$  & 96 & (i) &4.5  &3.3 &4.3 & 4.7  \\
&&& (ii) &18.2 & 16.8 & 27.2 & 45.8  \\
&&& (iii) &24.9 & 21.0 & 49.0 & 71.2  \\
&&& (iv) &32.6 & 30.6 & 75.9 & 61.0  \\
&&& (v) &3.2 &12.4  &4.8  & 6.8\\
B & $\mathbf{v}_X=(4,0.5,2,0.3)$ & 48 & (i) &8.2 &3.8  &7.0 & 6.1  \\
&&& (ii) & 43.3 & 45.8 & 83.4 & 71.4\\
&&& (iii) & 66.6 & 65.2 & 23.8 & 18.1\\
&&& (iv) & 26.5 & 22.5 & 5.8 & 3.7\\
&&& (v) &2.4 &6.0  &3.6  & 1.6\\
& $\mathbf{v}_Y=(2,0.5,1,0.3)$ & 96 & (i) &5.3  &4.5  &5.2& 4.9 \\
&&& (ii)  & 20.3 & 24.2 & 36.7 & 54.6 \\
&&& (iii) & 25.3 & 27.9 & 53.0 & 75.7\\
&&& (iv) & 33.1 & 32.6 & 78.1 & 60.1\\
&&& (v) &2.8 &8.4  &6.8  & 3.6\\
\midrule
\end{tabular*}
\end{table}

\section{Climate projections analysis}\label{data}

\begin{figure}

\includegraphics{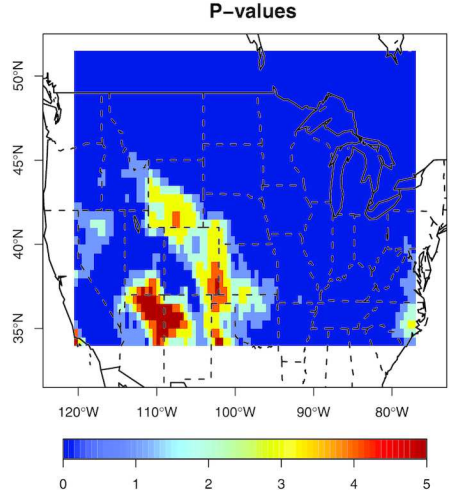}

\caption{$p$-values for testing the nullity of lag zero
cross-correlation between the station observations and model outputs
at each location. The numbers 0--5 denote the ranges of the
$p$-values, that is, 0 denotes [0.1, 1]; 1 denotes [0.05, 0.1]; 2
denotes [0.025, 0.05]; 3 denotes [0.01, 0.025]; 4 denotes [0.005,
0.01] and 5 denotes [0, 0.005].}\label{fig:pvalues}
\end{figure}

We apply the SN-based test to a gridded spatio-temporal temperature
dataset covering a subregion of North America. The dataset comes
from two separate sources: gridded observations generated from
interpolation of station records (HadCRU), and gridded simulations
generated by an AOGCM (NOAA GFDL CM2.1). Both datasets provide
monthly average temperature for the same 19-year period, 1980--1998.
Each surface is viewed as a two-dimensional functional datum. The
yearly average data have been recently analyzed in Zhang \textit{et al.}
\cite{r26}, where the goal is to detect a possible change point of the
bias between the station observations and model outputs. In this
paper, we analyze the monthly data from 1980 to 1998, which includes
228 functional images for each sequence. We focus on the second-order properties and aim to test the equality of the eigenvalues and
eigenfunctions of the station observations and model outputs. To
perform the analysis, we first remove the seasonal mean functions
from the two functional sequences. At each location, we have two
time series from the demeaned functional sequences. We apply the
SN-based test developed in Shao \cite{r24} to test
whether their
cross-correlation at lag zero is equal to zero. The $p$-values of
these tests are plotted in Figure~\ref{fig:pvalues}. The result
tends to suggest that the dependence between the station
observations and model outputs may not be negligible at certain
regions as the corresponding $p$-values are extremely small. The two
sample tests introduced in this paper are useful in this case
because they are robust to such dependence.

\begin{figure}

\includegraphics{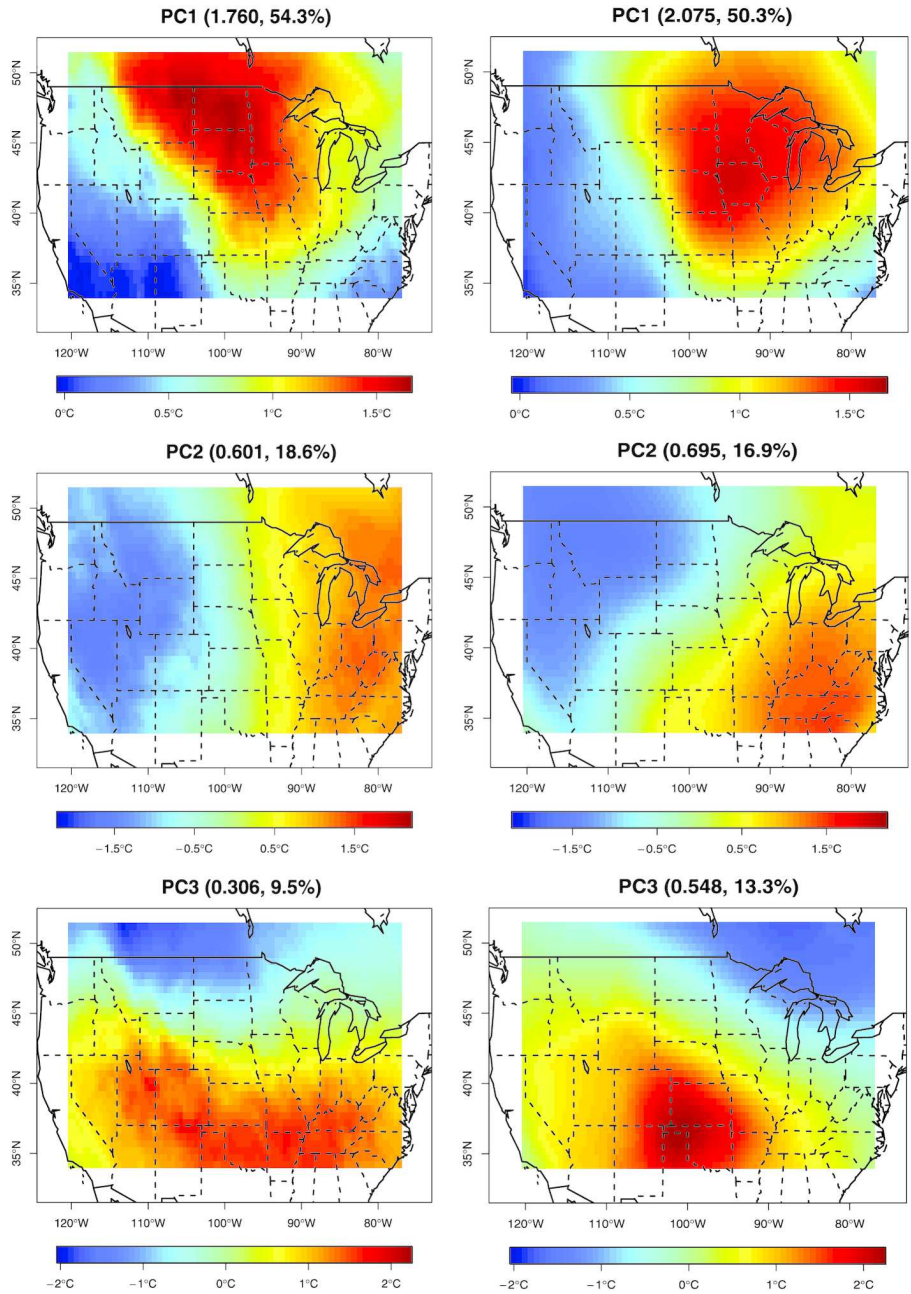}

\caption{The first three PCs of the station observations (left
panels) and model outputs (right panels), and the associated
eigenvalues and percentage of variations explained.}\label{fig:Pc}
\end{figure}

We perform FPCA on the demeaned sequences. Figure~\ref{fig:Pc} plots
the first three PC's of the station observations and model outputs.
We then apply the SN-based tests $G_{SN,N}^{(2)}(M)$ and
$G_{SN,N}^{(3)}(M_0)$ (with $p=3$) to the demeaned sequences, which
yields the results summarized in Table~\ref{tb:gcm}. It is seen from
the table that the first two eigenvalues of the station observations
and model outputs may be the same, at least statistical significance
is below the 10\% level, while there is a significant difference
between their third eigenvalue. The SN-based tests also suggest that
there are significant differences of the first and second PCs
between the station observations and model outputs as the
corresponding $p$-values are less than 5\% while the difference
between the third PCs is not significant at the 10\% level; compare
Figure~\ref{fig:Pc}. We\vspace*{-2pt} also tried the basis functions $\nu^*_j$ and
$\nu^{**}_j$ for $G_{SN,N}^{(3)}(M_0)$ (see (\ref{nu1}) and
(\ref{nu2})), which leads to the same conclusion. To sum up, our
results suggest that the second-order properties of the station
observations and model outputs may not be the same.

In climate projection studies, the use of numerical models outputs
has become quite common nowadays because of advances in computing
power and efficient numerical algorithms. As mentioned in Jun
\textit{et al.} \cite{r13}, ``Climate models are evaluated on how well they
simulate the
current mean climate state, how they can reproduce the observed
climate change over the last century, how well they simulate
specific processes, and how well they agree with proxy data for very
different time periods in the past.'' Furthermore, different
institutions produce different model outputs based on different
choices of parametrizations, model components, as well as initial
and boundary conditions. Thus there is a critical need to assess the
discrepancy/similarity between numerical model outputs and real
observations, as well as among various model outputs. The two sample
tests proposed here can be used towards this assessment at a
preliminary stage to get a quantitative idea of the difference,
followed by a detailed statistical characterization using
sophisticated spatio-temporal modeling techniques (see, e.g.,
Jun \textit{et al.} \cite{r13}). In particular, the observed
significance level for each
test can be used as a similarity index that measures the similarity
between numerical model outputs and real observations, and may be
used to rank model outputs. A detailed study along this line would
be interesting, but is beyond the scope of this article.

\begin{table}
\tablewidth=\textwidth
\tabcolsep=0pt
\caption{Comparison of the eigenvalues and eigenfunctions of the
covariance operators of the station observations and model outputs}\label{tb:gcm}
\begin{tabular*}{\textwidth}{@{\extracolsep{\fill}}lllll@{}}
\hline
$K$ & \multicolumn{1}{l}{$G_{SN,N}^{(2)}(M)$} & $p$-value & \multicolumn{1}{l}{$G_{SN,N}^{(3)}(M_0)$} & $p$-value \\
\hline
1    & \phantom{0}10.8  & $(0.1, 1)$ &  126.4 & $(0.025, 0.05)$\\
2    & \phantom{00}5.4  & $(0.1, 1)$ &  295.4   & $(0, 0.005)$\\
3    & 119.9 & $(0.005, 0.01)$  &  \phantom{0}34.2  & $(0.1, 1)$\\
--   & 326.2 & $(0.005, 0.01)$ & 318.0 & $(0.005, 0.01)$\\
\hline
\end{tabular*}
\legend{Note: The first three rows show the results for testing individual
eigencomponent, and the last row shows the results for testing the
first three eigenvalues and eigenfunctions jointly.}
\end{table}

\section*{Acknowledgements}

We are grateful to an associate editor and a referee for helpful
comments, which led to
substantial improvements. Shao's research is supported in part by National
Science Foundation Grant DMS-11-04545.

\begin{supplement}
\stitle{Supplement to ``Two sample inference for the second-order property
of temporally dependent functional data''}
\slink[doi]{10.3150/13-BEJ592SUPP} 
\sdatatype{.pdf}
\sfilename{BEJ592\_supp.pdf}
\sdescription{This supplement contains proofs of the results in
Section~\ref{theory}.}
\end{supplement}

%

\printhistory

\end{document}